\documentclass[11pt]{amsart}

\usepackage[all]{xypic}

\usepackage{latexsym}

\usepackage{amssymb}
\usepackage{amsfonts}
\usepackage{amscd}
\usepackage{amsmath,amsthm}

\newtheorem{lemma}{Lemma}[section]

\newtheorem{lem}[lemma]{Lemma}
\newtheorem{prop}[lemma]{Proposition}
\newtheorem{thm}[lemma]{Theorem}
\newtheorem{cor}[lemma]{Corollary}

{

}

\theoremstyle{definition}

\theoremstyle{remark}

\numberwithin{equation}{section}

\newenvironment{pf}{\noindent{\bf Proof.}}{\hfill $\square$\medskip}

\def\CC{{\mathbb C}}

\def\PP{{\mathbb P}}

\def\XX{{\mathbb X}}
\def\ZZ{{\mathbb Z}}

\def\0ol{{\bar 0}}
\def\1ol{{\bar 1}}
\def\2ol{{\bar 2}}
\def\ol2{{\bar 2}}
\def\3ol{{\bar 3}}
\def\4ol{{\bar 4}}
\def\5ol{{\bar 5}}
\def\6ol{{\bar 6}}
\def\7ol{{\bar 7}}
\def\8ol{{\bar 8}}
\def\9ol{{\bar 9}}

\def\bold0{{\bf 0}}
\def\bold1{{\bf 1}}
\def\bold2{{\bf 2}} 
\def\bold3{{\bf  3}}
\def\bold4{{\bf 4}}
\def\bold5{{\bf 5}}
\def\bold6{{\bf 6}}
\def\bold7{{\bf 7}}
\def\bold8{{\bf 8}}
\def\bold9{{\bf 9}}

\def\xul{{\underline{x}}}

\def\P2Skly{\PP^2_{Skly}}

\def\coker{\operatorname {coker}}

\def\End{\operatorname {End}}
\def\Ext{\operatorname {Ext}}

\def\gr{\operatorname {gr}}

\def\Hom{\operatorname {Hom}}

\def\Tor{\operatorname {Tor}}

\def\dim{\operatorname{dim}}

\def\End{\operatorname{End}}

\def\Ext{\operatorname{Ext}}
\def\Fdim{{\sf Fdim}}
\def\fdim{{\sf fdim}}

\def\Gr{{\sf Gr}}

\def\gr{{\sf gr}}

\def\Hom{\operatorname{Hom}}

\def\id{\operatorname{id}}

\def\liminj{\varinjlim}

\def\mod{{\sf mod}}
\def\Mod{{\sf Mod}}

\def\Projnc{\operatorname{Proj}_{nc}}

\def\QGr{\operatorname{\sf QGr}}
\def\qgr{\operatorname{\sf qgr}}

\def\rank{\operatorname{rank}}

\def\uExt{\operatorname{\underline{Ext}}}

\def\uHom{\operatorname{\underline{Hom}}}

\def\ul1{\operatorname{\underline{1}}}

\def\d{\downarrow}

\def\a{\alpha}
\def\b{\beta}

\def\d{\delta}

\def\sA{{\sf A}}

\def\sD{{\sf D}}

\def\sT{{\sf T}}

\def\cal{\mathcal}

\def\cE{{\cal E}}
\def\cF{{\cal F}}
\def\cG{{\cal G}}

\def\cM{{\cal M}}

\def\cO{{\cal O}}

\def\coh{{\sf coh}}

\def\Qcoh{{\sf Qcoh}}

\def\dirlim{\mathop{\vtop{\baselineskip -100pt\lineskip -1pt\lineskiplimit 0pt
\setbox0\hbox{lim}\copy0\hbox to \wd0{\rightarrowfill}}}\limits}
\def\invlim{\mathop{\vtop{\baselineskip -100pt\lineskip -1pt\lineskiplimit 0pt
\setbox0\hbox{lim}\copy0\hbox to \wd0{\leftarrowfill}}}\limits}

\def\I11{{1 \kern -0.8pt \! \mbox{l}}}
\def\mumu{{\mu\kern-4.2pt\mu}}
\def\bfmu{{\mu\kern-4.2pt\mu}}
\def\2slash{\backslash \! \backslash}

\def\boxtimes{\setbox0\hbox{$\Box$}\copy0\kern-\wd0\hbox{$\times$}}

\date{}                                           

\begin{document}

\title[The free algebra as a homogeneous coordinate ring]{The non-commutative scheme having a free algebra as a homogeneous coordinate ring}

\author{S. Paul Smith}

\address{ Department of Mathematics, Box 354350, Univ.
Washington, Seattle, WA 98195}

\email{smith@math.washington.edu}

\thanks{This work was partially supported by NSF grant 0602347. }


\begin{abstract}
Let $k$ be a field and $TV$ the tensor algebra on a finite-dimensional $k$-vector space $V$. 
This paper proves that the quotient category $\QGr(TV) := \Gr(TV)/\Fdim$ of graded $TV$-modules modulo those that are unions of finite dimensional modules is equivalent to the category of modules over the direct limit of matrix algebras, $\liminj_r M_n(k)^{\otimes r}$. Non-commutative algebraic geometry associates to a graded algebra $A$ a ``non-commutative scheme'' $\Projnc A$ that is defined implicitly by declaring that the 
category of ``quasi-coherent sheaves'' on $\Projnc A$ is $\QGr A$. When $A$ is coherent and $\gr A$ its category of finitely presented graded modules, $\qgr A:= \gr A/\fdim$ is viewed as the category of ``coherent sheaves'' on $\Projnc A$. We show that when $\dim V \ge 2$, $\qgr(TV)$ has no indecomposable 
objects, no noetherian objects, and no simple objects.  Moreover, every short exact sequence in $\qgr(TV)$ splits.  

We also prove $\QGr(TV) \equiv \Gr L$ where $L$ is the Leavitt algebra on $2\dim V$ generators that 
embeds as a dense subalgebra of the Cuntz algebra $\cO_{\dim V}$.

\end{abstract}

\maketitle

\pagenumbering{arabic}


\setcounter{section}{0}

\section{Introduction}

\subsection{}
\label{sect.1.1}
Let $n$ be a positive integer. 

Throughout this paper $k$ is a field and 
$$
R:=k\langle x_0, x_1,\ldots,x_n\rangle
$$
is the free algebra on $n+1$ variables with $\ZZ$-grading given by declaring that $\deg x_i=1$ for all $i$.
This paper concerns the categories of coherent and quasi-coherent ``sheaves''
on the ``non-commutative scheme''
$$
\XX^n:=\Projnc k\langle x_0,x_1,\ldots,x_n\rangle
$$
with ``homogeneous coordinate ring'' $R$.

The ``scheme'' $\Projnc R$ is an imaginary object: there is no underlying topological 
space endowed with a sheaf of rings. Rather one declares that the category of ``quasi-coherent sheaves'' on $\Projnc R$ is the quotient category
$$
\Qcoh(\Projnc R):=\QGr R:=\frac{\Gr R}{\Fdim  R}
$$
where $\Gr R$ is the category of $\ZZ$-graded left $R$-modules with degree-preserving homomorphisms
 and $\Fdim R$ its full subcategory of direct limits of finite-dimensional modules. The imaginary space 
 $\Projnc R$  manifests itself through the category $\Qcoh(\Projnc R)$.

The category $\QGr R$ and its full subcategory of finitely presented objects,  $\qgr R$,  is
the focus of this paper.\footnote{An object $\cM$ in an additive category $\sA$ is finitely presented if $\Hom_\sA(\cM,-)$ commutes with direct limits; is finitely generated if whenever $\cM = \sum \cM_i$ for some directed family of subobjects $\cM_i$ there is an index $j$ such that $\cM=\cM_j$; is coherent if it is finitely presented and all its finitely generated subobjects are finitely presented.} We think of $\qgr R$ as the category of ``coherent sheaves'' on $\Projnc R$.

\subsection{}

The notations $\QGr R$ and $\Qcoh \XX^n$ are interchangeable as are the notations $\qgr R$ and $\coh \XX^n$. The reader may adopt either so as to reinforce either an algebraic or 
geometric perspective. 

\subsection{}

We write 
$$
\pi^*:\Gr R \to \QGr R
$$
for the quotient functor and define $\cO:=\pi^* R$. We call $\cO$ a {\sf structure sheaf} for $\Projnc R$. 
 
 For the commutative polynomial ring $k[x_0,\ldots,x_n]$, 
 $\QGr k[x_0,\ldots,x_n]$ is equivalent to  $\Qcoh \PP^n$ and $\pi^*$
``is'' the functor usually denoted $M \mapsto \widetilde{M}$ in algebraic geometry texts and $\pi^*(k[x_0,\ldots,x_n])$ {\it is} the   structure sheaf $\cO_{\PP^n}$.

\subsection{}

 The {\sf twist functor} on $\Gr R$, denoted
$M \rightsquigarrow M(n)$, 
 is defined by $M(n)_i:=M_{n+i}$ with the same action of $R$.  The subcategory $\Fdim R$ is stable under
 twisting so there is an induced functor on $\QGr R$ that we denote by 
 $\cF \rightsquigarrow  \cF(n)$ and call the {\sf Serre twist}.

\subsection{The main results}

We define the direct limit algebra
$$
S:= \liminj_i M_{n+1}(k)^{\otimes i}
$$
where the maps in the directed system are $a_1 \otimes \cdots \otimes a_i \mapsto 1 \otimes a_1 \otimes \cdots \otimes a_i  $. 
The ring $S$ is coherent so finitely presented $S$-modules form an abelian category.
 
  \begin{thm}
  \label{thm.3}
There is an equivalence of categories
$$
\Hom_{\QGr R}(\cO,-):\QGr R \equiv \Mod_r S,
$$
the category of right $S$-modules.
The equivalence sends $\cO$ to $S$, i.e., $S=\Hom_{\QGr R}(\cO,\cO)$. Furthermore,
the equivalence restricts to an equivalence
$$
\qgr R \equiv \mod_r S,
$$
the category of finitely presented right $S$-modules.
\end{thm}

The ring $S$ is anti-isomorphic to itself so it doesn't really matter whether we choose to work with left or right $S$-modules.

The key to Theorem \ref{thm.3} is the following preliminary result.

\begin{thm}
  $\cO$ is a finitely generated, projective, generator in $\Qcoh \XX^n$.
  \end{thm}
  
As the next result emphasizes, the categories $\coh \XX^n$ and $\Qcoh \XX^n$ are unlike 
the categories of  coherent and quasi-coherent sheaves over quasi-projective schemes.

\begin{thm}
Suppose $n \ge 1$. 
\begin{enumerate}
  \item 
There are no indecomposable objects in $\coh \XX^n$, hence no simple objects, and therefore 
no noetherian objects other than 0.
  \item 
   Every short exact sequence in $\coh \XX^n$ splits. 
  \item 
  Every object in $\coh \XX^n$ is isomorphic to a finite direct sum of various $\cO(i)$s with 
  finite multiplicities.
  \item{}
  If $\cF,\cG \in \coh \XX^n$ are non-zero, then $\dim_k \Hom_{\XX^n}(\cF,\cG) = \infty$.
  \item{}
  The Grothendieck group of the abelian category $\coh \XX^n$ 
is isomorphic to $\ZZ\big[\frac{1}{n+1}\big]$ as an additive group. 
\end{enumerate}
\end{thm}

The reason $\coh \XX^n$ and $\Qcoh \XX^n$ behave so differently from categories of (quasi-)coherent 
sheaves on quasi-projective schemes is that $\cO(-1)$ is a non-trivial direct summand of $\cO$. Since $R_{\ge 1}$ is isomorphic to $R(-1)^{\oplus(n+1)}$, 
$$
\cO \cong \cO(-1)^{\oplus(n+1)}.
$$ 
This behavior also occurs for the graded algebras $k\langle x,y\rangle/(y^{r+1})$  in \cite{Sm.y^2=0}. 

Section 2 of \cite{Sm.y^2=0} establishes some general results that are applied in the present paper to the free algebra: almost everything in the present paper is a rather simple consequence of those more general results.

 \subsection{Connection to Leavitt algebras and Cuntz algebras}
 Let $L$ be the Leavitt algebra generated by the entries in the row vectors $\xul=(x_0,\ldots,x_n)$
 and $\xul^*=(x_0^*,\ldots,x_n^*)$ subject to the relations 
 $$
 \xul^*\xul^{\sT}=1
 \qquad \hbox{and} \qquad
 \xul^{\sT}\xul^* =I_{n+1},
 $$
 the $(n+1) \times (n+1)$ identity matrix. Give $L$ a $\ZZ$-grading by $\deg x_i=1$ and $\deg x_i^*=-1$. 
 The connection between $\QGr R$ and $L$ is made in the following result which is proved in section 4.
 
 \begin{thm}
 \label{thm.GrL}
 The algebra $L$ is strongly graded, $L_0 \cong S$, and there is an equivalence of categories
 $$
 \QGr R \equiv \Gr L \equiv \Mod L_0.
 $$
 \end{thm}
 
 The proof makes use of the following facts: $L$ is the universal localization of $R$ that inverts the  homomorphism $R(-1)^{n+1} \to R$ whose cokernel is $R/R_{\ge 1}$; $L$ is flat as a right $R$-module; if $M \in \Gr R$ and $L \otimes_R M=0$, then $M \in \Fdim R$. The first of these facts is well-known; the second is proved in \cite{A}; a version of the third for modules in $\gr R$ is proved in \cite{A}.

  In \cite{Cuntz}, Cuntz defined a class of C$^*$-algebras $\cO_{n+1}$, $n \ge 1$, generated by the
  ``same'' elements and relations for $L$. The algebra $L$ is a dense subalgebra of $\cO_{n+1}$.

\medskip
{\bf Acknowledgements.}
The results in this note grew out of conversations at the workshop  {\it Test problems for the theory of finite dimensional algebras} held at the Banff International Research Station in September 2010.  I  thank the organizers of the conference for arranging a stimulating meeting. Particular thanks to Helmut Lenzing for his questions and interest in this work. I am grateful to Ken Goodearl for 
answering questions about von Neumann regular rings, AF algebras, and their Grothendieck groups. 

The results in section 4 were stimulated by lectures given by Pere Ara at the Latin-American School of Mathematics held in Cordoba, Argentina, in May 2011, and conversations with Pere Ara, Ken Goodearl, and Michel Van den Bergh. I thank them all.

\section{Coherent sheaves on $\Projnc TV$}
\label{sect.proj.TV}

Let $V$ be a $k$-vector space of dimension $d=n+1$ and 
$$
R:=TV=k\langle x_0,\ldots,x_n\rangle
$$
be the tensor algebra on $V$ with $\ZZ$-grading  $R_i:= V^{\otimes i}$.

 \subsection{Graded modules over the free algebra}
A ring is {\sf left coherent} if all its finitely generated left ideals are finitely presented. Since $TV$ is anti-isomorphic to itself we can dispense with the adjectives {\it left} and {\it right} when discussing properties
of $TV$ like coherence. 

The important property for us, indeed an equivalent characterization of coherence, is that the category of finitely presented left modules over a left coherent ring is abelian.

 The following facts are well-known.
 
 \begin{prop}
 Let $V$ be a vector space of finite dimension $d \ge 1$.
 \begin{enumerate}
  \item
  Every left ideal in $TV$ is free.
  \item 
   $TV$ is coherent and has global dimension one. 
  \item 
   Every finitely generated projective $R$-module is free.
   \item{}
   If $d \ge 2$, $TV$ has exponential growth:  $H(R,t)=(1-dt)^{-1}$. 
      \item{}
   If $d \ge 2$, $TV$ is not noetherian.
   \item{}
    $R_{\ge i}$ is isomorphic to $R(-i)^{d^i}$, the free $R$-module of rank 
     $d^i$ with basis in degree $i$. 
\end{enumerate}
 \end{prop}
 
Suppose $N$ is a graded $R$-module. There is an exact sequence $L \to M \to N \to 0$ in $\Mod R$
with $L$ and $M$ finitely generated  if and only if there is an exact sequence  
$L' \to M' \to N \to 0$ in $\Gr R$ with $L'$ and $M'$ generated by a finite number of homogeneous elements. 

We define
$$
\gr R:=\hbox{the category of finitely presented graded $R$-modules}.
$$
 This is an abelian category because $R$ is coherent.

 \begin{prop}
 \label{prop.gr.coh}
 Let $R=TV$ and $M \in \QGr R$.
 \begin{enumerate}
  \item 
 $M$ is graded-coherent if and only if for all $i \gg 0$, 
 $$
 M_{\ge i} \cong R(-i)^{t_i}
 $$
 for some integer $t_i$ depending on $i$. 
  \item 
  If $0 \to L \to M \to N \to 0$ is an exact sequence in $\gr R$, then 
   $0 \to L_{\ge i} \to M_{\ge i} \to N_{\ge i} \to 0$ splits for $i \gg 0$.
    \item 
  If $M \in \gr R$, then $M$ has a largest finite dimensional graded submodule.
\end{enumerate}
 \end{prop}
 \begin{pf}
 (1)
 Suppose $M$ is finitely presented. Then there is an exact sequence $0 \to F' \to F \to M \to 0$ in $\gr R$ with $F$ and $F'$ finitely generated graded free $R$-modules. Since $F'$, $F$, and $M$, are finitely generated graded modules, 
 $F'_{\ge i}$, $F_{\ge i}$, and $M_{\ge i}$, are generated as $R$-modules by $F'_i$, $F_i$, and $M_i$,
 respectively for all sufficiently large $i$.  But $R_{\ge m} \cong R(-m)^{d^m}$ so
 $$
 \big(0 \to F' \to F \to M \to 0\big)_{\ge i} \; = \; \big( 0 \to R(-i)^r \to R(-i)^{s} \to M_{\ge i} \to 0 \big)
 $$
 for $i \gg 0$. 
Every degree-zero homomorphism $R(-i)^r \to R(-i)^{s}$ splits so $M_{\ge i} \cong R(-i)^{s-r}$.
 
 The converse is trivial.
 
 (2) 
 By (1), $N_{\ge i}$ is free for $i \gg 0$, hence the splitting.
 
 (3)
 Since $R$ is a domain its only finite dimensional submodule is zero.
 It now follows from (1) that the only finite dimensional submodule of $M_{\ge i}$
 is the zero submodule for $i \gg 0$. There is therefore an integer $n$ such that every finite dimensional
 submodule of $M$ is contained in $\sum_{j=-i}^i M_j$. But $M$ is finitely generated so $\sum_{j=-i}^i M_j$
 has finite dimension. Hence the sum of all finite dimensional submodules of $M$ has finite dimension, and that sum is therefore the largest finite dimensional graded submodule of $M$. 
  \end{pf}

   \subsection{}
We write $\fdim R$ for $\Fdim R \cap \gr R$. Thus $\fdim R$ is the full subcategory of $\Gr R$ consisting of the finite dimensional submodules. We define the category of ``coherent sheaves'' on $\XX^n$ by
$$
\coh \XX^n=\qgr R:= \frac{\gr R}{\fdim R}.
$$
Since $\Fdim$ satisfies condition (2) of \cite[Prop. A.4, p. 113]{HK} with respect to the Serre subcategory $\fdim$ of $\gr R$, $\Fdim$ is localizing of finite type which then allows us to apply  \cite[Prop. A.5, p. 113]{HK}
and so conclude that
$\coh \XX^n$ consists of finitely presented objects in $\Qcoh \XX^n$ and every object in  $\Qcoh \XX^n$ is a direct limit of objects in $\coh \XX^n$. 
 
\begin{prop}
\label{prop.ss}
Every short exact sequence in $\qgr R$ splits.
\end{prop}
\begin{pf}
By \cite[Cor. 1, p. 368]{Gab}, every short exact sequence in $\qgr R$  is of the form
$$
0 \longrightarrow  \pi^* L \stackrel{\pi^*f}{\longrightarrow}  \pi^* M \stackrel{\pi^*g}{\longrightarrow}  \pi^* N \longrightarrow 0
$$
where $0 \longrightarrow  L \stackrel{f}{\longrightarrow}  M \stackrel{g}{\longrightarrow}   N \longrightarrow 0$ 
is an exact sequence in $\gr R$. But 
$$
\pi^* \big(0 \to L \to  M \to  N \to 0\big) \, \cong \,  \pi^* \big(0 \to L_{\ge i} \to  M_{\ge i} \to  N_{\ge i}\to 0\big)
$$
for all $n$ and $N_{\ge i}$ is free for $i \gg 0$ so the sequence $0 \to L_{\ge i} \to  M_{\ge i} \to  N_{\ge i}\to 0$ splits for $i \gg 0$. Applying $\pi^*$ to a split exact sequence yields a split exact sequence, whence the result.
\end{pf}

 \begin{cor}
 \label{cor.classfcn1}
Let $\cF \in \qgr R$. If $i \gg 0$ there is an integer $r$, depending on $i$,  such that
$$
\cF \cong \cO(i)^r.
$$
 \end{cor}
 \begin{pf}
 There is an $M \in \gr R$ such that $\cF = \pi^*M$. But 
 $\dim_k M/M_{\ge i} < \infty$ so $\pi^*M \cong \pi^* M_{\ge i}$ for all $i \gg 0$. Now apply Proposition
 \ref{prop.gr.coh}(1).
 \end{pf}

\begin{cor}
\label{cor.classfcn2}
Every object in $\qgr R$ is  injective and projective. Furthermore, $\cO$ is projective as an object in $\QGr R$.
\end{cor}
\begin{pf}
It follows from  Proposition \ref{prop.ss} that every object in $\qgr TV$ is injective and projective in $\qgr TV$. 
Since every left ideal of $R$ is free it is a tautology that  every graded left ideal that has finite codimension in $R$ contains a {\it free} graded left ideal that has finite codimension in $R$. Hence 
\cite[Prop. 2.7]{Sm.y^2=0} applies and tells us that $\cO$ is projective  in $\QGr R$.
\end{pf}

 \subsection{}
Since $R_{\ge i} \cong R(-i)^{\oplus (\dim V)^i}$ has finite codimension in $R$, 
\cite[Prop. 2.7 and Thm. 2.8]{Sm.y^2=0} may be applied to yield the next two results.   
 
See Footnote 1 for the definition of finitely generated, finitely presented, and coherent, objects in an abelian category. 

\begin{prop}
\label{prop.fpres}
Let $A$ be a left graded-coherent ring and write $\cO$ for the image of $A$ in $\QGr A$. Then 
\begin{enumerate}
  \item 
  $\QGr A$ is a locally coherent category; 
  \item 
  the  full subcategory of $\QGr A$ consisting of the finitely presented objects is 
equivalent to $\qgr A$;
  \item 
$\cO$ is coherent: i.e.,  $\Hom_{\QGr A}(\cO,-)$ commutes with direct limits;
\item{}
 $\cO$ is finitely generated.
\end{enumerate}
\end{prop}
 
\begin{thm}
\label{thm.progenor}
$\phantom{xx}$
\begin{enumerate}
  \item 
  $\cO$ is a progenerator in $\QGr R$.
  \item 
The functor $\Hom(\cO,-)$ is an equivalence from the category $\QGr R$  to the category of right modules over the endomorphism ring $\End_{\QGr R} \cO$.
\item{}
$\End_{\QGr R} \cO \cong \liminj \End_{\Gr R}(R_{\ge i})$, the direct limit of the directed system
\begin{equation}
\label{eq.dsys.End}
  \UseComputerModernTips
\xymatrix{
\cdots \ar[r] & \End_{\Gr R}(R_{\ge i})   \ar[r]^{\theta_i} & \End_{\Gr R}(R_{\ge i+1}) \ar[r] & \cdots
}
\end{equation}
of $k$-algebras in which $\theta_i(f) = f|_{R_{\ge i+1}}$. 
\end{enumerate}
\end{thm}

We will determine   this direct limit  in section \ref{sect.The.dlim.S}.
 
 \subsection{}
 
  \begin{lem}
  \label{lem.decomposing}
If $\dim V=d \ge 1$, then
$$
\cO \cong \cO(-1)^{\oplus d} \cong \cO(-2)^{\oplus d^2} \cong \cdots
$$ 
 \end{lem}
 \begin{pf}
From the exact sequence
$
0 \to R_{\ge 1} \to R \to k \to 0
$
we see that $\pi^* R \cong \pi^* R_{\ge 1}$. But $R_{\ge 1} \cong R \otimes_k V \cong R(-1)^d$  so
$\pi^* R_{\ge 1} \cong \cO(-1)^d$. Hence $\cO=\pi^*R \cong \cO(-1)^d$. The result now follows by induction.
 \end{pf}
 
Lemma \ref{lem.decomposing} implies that for every $\cF \in \qgr R$ and every $r \ge 1$, 
 $$
 \cF \cong \cF(-r)^{\oplus d^r}.
 $$

 \begin{cor}
 \label{cor.matrices.endoms}
 Let $d = \dim V \ge 1$. If $\cF \in \qgr R$ is non-zero, then for all integers $r \ge 0$ there is an injective ring homomorphism
$$
M_{d^r}(k) \to \End \cF.
$$
 \end{cor}

\begin{cor}
Suppose $\dim V \ge 2$. Then
\begin{enumerate}
  \item 
  the only noetherian object in $\qgr R$ is the zero object; 
  \item 
there are no simple objects in $\qgr R$;
\item{}
there are no indecomposable objects in $\qgr R$;
\item{}
$\dim_k \Hom_{\QGr R} (\cF,\cG)=\infty$ for all non-zero objects $\cF$ and $\cG$ in $\qgr R$.
\end{enumerate}
\end{cor}
\begin{pf}
(2)
By Schur's lemma the endomorphism ring of a simple object is a division algebra but 
Corollary \ref{cor.matrices.endoms} says that endomorphism rings of objects in $\coh \XX^{n}$ are 
never division rings when $\dim V \ge 2$. Hence  $\coh \XX^{n}$ has no simple objects. Part (1) follows immediately because a non-zero noetherian object has at least one maximal subobject and hence 
a simple quotient.

(3)
See the remark after Lemma \ref{lem.decomposing}.

(4) 
Since $\Hom_{\XX^{n}}(\cF,\cG)$ is an $\End \cG$-$\End\cF$-bimodule it is a module over the matrix algebra
$M_{d^r}(k)$ for all $r \ne 0$. It therefore suffices to show that $\Hom_{\XX^{n}}(\cF,\cG)$ is non-zero. 
By Lemma \ref{lem.decomposing} there is a monic map $\cO(-1) \to \cO$ and an epic map $\cO \to \cO(-1)$.
Because the twist $(1)$ is an auto-equivalence it follows (using compositions of twists of the monic and epic maps just mentioned) that there is a monic map $\cO(i) \to \cO(j)$ whenever $i \le j$ and an 
epic map $\cO(j) \to \cO(i)$ whenever $j \ge i$.  It now follows from Corollary \ref{cor.classfcn1}
that $\Hom_{\XX^{n}}(\cF,\cG) \ne 0$. 
\end{pf}

\section{A direct limit of matrix algebras}
\label{sect.The.dlim.S}

We will now show that the endomorphism ring of $\cO$ is isomorphic to the ring $S$ we are about to define.

Let $S_i:=M_{n+1}(k)^{\otimes i}$ and define 
\begin{equation}
\label{eq.d.sys}
\theta_i: S_i \to S_{i+1}=M_{n+1}(k)\otimes S_i \qquad \hbox{by} \qquad  \theta_i(a)=1 \otimes a.
\end{equation}
The homomorphisms $\theta_i$ determine a directed system and we define
$$
S:=\liminj_i S_i.
$$

We  write  $\Mod S$ for the category of right $S$-modules and $\mod S$ for its full 
category of finitely presented $S$-modules.

\begin{thm}
\label{thm.equiv.cats}
If $S$ is the ring above, then $\End_{\XX^{n}}\cO \cong S$,  
the functor $\Hom_{\XX^{n}}(\cO,-)$ is an equivalence of categories
$$
\Qcoh \XX^{n} \equiv \Mod S
$$
sending $\cO$ to $S_S$, and  $\Hom_{\XX^{n}}(\cO,-)$ restricts to an equivalence
$$
\coh \XX^{n} \equiv \mod S.
$$
\end{thm} 
\begin{pf}
Write $R=TV$ as in section \ref{sect.proj.TV}. 
By the definition of  morphisms in a quotient category,
\begin{equation}
\label{.d.sys1}
\End_{\XX^{n} }\cO =\Hom_{\QGr R}(\pi^* R , \pi^* R) = \liminj \Hom_{\Gr R}(R',R/R'')
\end{equation}
where $R'$ runs over all graded left ideals in $R$ such that $\dim_k(R/R') < \infty$ and 
$R''$ runs over all graded left ideals in $R$ such that $\dim_k R'' < \infty$. 

By Theorem \ref{thm.progenor} (see \cite[Sect. 2]{Sm.y^2=0} for a fuller explanation), this reduces to
$$
\End_{\QGr R} \cO = \liminj_{i} \Hom_{\Gr R}(R_{\ge i},R_{\ge i}).
$$

As a left $R$-module, $R_{\ge i} \cong R \otimes_k V^{\otimes i}$ where $V$ is placed in degree $1$. 
There is a commutative diagram
$$
\UseComputerModernTips
\xymatrix{
\Hom_k (V^{\otimes i} , V^{\otimes i})\ar[d]_{\rho_i}  \ar[rr]^{\psi_i}  && \Hom_k (V^{\otimes i+1} , V^{\otimes i+1})  \ar[d]^{\rho_{i+1}}
\\
 \Hom_{\Gr R}(R_{\ge i},R_{\ge i}) \ar[rr]_{\theta_i} &&   \Hom_{\Gr R}(R_{\ge i+1},R_{\ge i+1})
 }
 $$
 in which 
 $$
 \rho_i(f)(r \otimes v)=rf(v) \qquad \hbox{ for all $r \in R$ and $v \in V^{\otimes i}$}
$$
and
 $$
 \psi_i(f) (a_0 \otimes a_1 \otimes \cdots \otimes a_i) = a_0 \otimes f(a_1 \otimes \cdots \otimes a_i).
 $$
But $\rho_i$ is an isomorphism so  
$$
\End_{\XX^{n}} \cO \cong  \liminj_{i}  \End_k  V^{\otimes i}  \cong \liminj_{i} M_{n+1}(k)^{\otimes i}.
$$
The proof is complete.
\end{pf}

The following properties of $S$ are either obvious or well-known (see \cite{KG2}). 

\begin{prop}
\label{prop.S}
$\phantom{xx}$
\begin{enumerate}
\item{}
$S$ is a simple ring;
\item{}
$S$ has no finite dimensional modules other than 0.
\item{}
Every finitely generated left ideal of $S$ is generated by an idempotent.
\item 
$S$ is a von Neumann regular ring.
\item 
$S$ is left and right coherent.
  \item 
  Every left $S$-module is flat.
    \item 
  Every finitely generated left $S$-module is projective.
  \item{}
  $K_0(S) \cong  \ZZ\big[\frac{1}{n+1}\big]$ via an isomorphism sending $[S]$ to 1.
\end{enumerate}
\end{prop}
\begin{pf}
(8)
Since $K_0(-)$ commutes with direct limits and $K_0(S_i)\cong \ZZ$ with $[S_i]=n+1$
under the isomorphism, $K_0(S)$ is the direct limit of the directed system
$$
\cdots \longrightarrow \ZZ \stackrel{n+1}{\longrightarrow}  \ZZ \stackrel{n+1}{\longrightarrow} \cdots.
$$
This direct limit is obviously isomorphic to $\ZZ[\frac{1}{n+1}]$. 
\end{pf}

It seems worthwhile to determine the Grothendieck group of $\qgr R$ independently of the equivalence of categories, i.e., without appealing to part (8) of Proposition \ref{prop.S}.   The next result does this.

\begin{prop}
\label{prop.K0}
As an additive group, the Grothendieck group of $\coh\XX^{n}$ is isomorphic to $\ZZ[\frac{1}{n+1}]$ with
$[\cO(i)] \longleftrightarrow (n+1)^i$.
\end{prop}
\begin{pf} 
Since $\fdim R$ is a Serre subcategory of $\gr R$ there is an exact sequence
$$
K_0(\fdim R) \to K_0(\gr R) \to K_0(\coh \XX^{n}) \to 0.
$$
It is clear that $K_0(\gr R) \cong \ZZ[t^{\pm 1}]$ via $R(-i) \leftrightarrow t^i$.
From the exact sequence $0 \to R(-1)^{n+1} \to R \to k \to 0$ we obtain $[k]=[R]-(n+1)[R(-1)]=1-(n+1)t$. 
Modules in $\fdim R$ are finite dimensional so have compositions series in which the composition factors are 
of the form $k(i)$ for various integers $i$. Hence, by d\'evissage,  $K_0(\fdim R) \cong K_0(\gr k)$ which is isomorphic as an additive group to $\ZZ[t^{\pm 1}]$ with $k(i) \longleftrightarrow t^{-i}$. The image of 
the map $K_0(\fdim R) \to K_0(\gr R)$ is therefore the ideal in $\ZZ[t^{\pm 1}]$ generated by the image of $[k]$.
Therefore 
$$
K_0(\coh \XX^{n}) \cong \frac{\ZZ[t]}{(1-(n+1)t)} \cong \ZZ\Big[\frac{1}{n+1}\Big].
$$
This completes the proof.
\end{pf}

\begin{cor}
If $m \ne n$, then $\XX^m \not\cong  \XX^n$.
\end{cor}

It is reasonable to define $\rank \cO(i) = (n+1)^i$.

\section{The relation to the Leavitt algebra $L(1,n+1)$ and the Cuntz algebra $\cO_{n+1}$}

As before $R$ is  the free algebra $k\langle x_0,\ldots,x_n\rangle$.

The main result in this section is that $\QGr R \equiv \Gr L \equiv  \Mod L_0$ where $L=L(1,n+1)$ is the 
finitely generated $\ZZ$-graded algebra defined below.  It is well-known that $L_0$ is the algebra we have called $S$ in the earlier part of this paper. 

\subsection{The Leavitt algebra}
The {\sf Leavitt algebra} $L=L(1,n+1)$, first defined in \cite{L},  is the $k$-algebra 
generated by  elements $x_0,\ldots,x_n,x_0^*,\ldots,x_n^*$ subject to the relations
\begin{equation}
\label{L.relns}
x_ix_i^*=1=x_0^*x_0+\cdots + x_n^*x_n
\qquad \hbox{and} \qquad x_ix_j^*=0 \; \; \hbox{if $i \ne j$}.
\end{equation}
A more meaningful definition is that $L$ is the universal localization \cite[Sect. 7.2]{Cohn-FRR}
of $R$ inverting the injective homomorphism 
$$
\iota:R^{n+1} \to R, \qquad (r_0,\ldots,r_n) \mapsto r_0 x_0+\cdots +r_nx_n.
$$
Since $\iota$ is right multiplication by $(x_0,\ldots,x_n)^\sT$, the formal inverse of $\iota$ is right multiplication by $(x_0^*,\ldots,x_n^*)$ where  
$$
(x_0,\ldots,x_n)^\sT(x_0^*,\ldots,x_n^*)=I_{n+1}
\quad \hbox{and} \quad (x_0^*,\ldots,x_n^*)(x_0,\ldots,x_n)^\sT =1
$$
and $I_{n+1}$ is the identity matrix.

\subsection{The Cuntz algebra}
The Cuntz algebra $\cO_{n+1}$ \cite{Cuntz} is the universal $C^*$-algebra generated by elements 
$x_0,\ldots,x_n$ subject to the relations (\ref{L.relns}). It is well-known that $L$ 
embeds in $\cO_{n+1}$ as a dense subalgebra. Much of the work in Cuntz's paper  \cite{Cuntz}  
involves purely algebraic calculations carried out inside $L$. 

\subsection{}
\label{sect.tors1}
One anticipates a relation between $\Gr L$ and $\QGr R$ because the fact that $\id_L \otimes \iota$ is an isomorphism implies that 
$$
0= L \otimes_R \coker (\iota) = L \otimes_R(R/R_{\ge 1}).
$$
It follows that $L \otimes_R -$ kills all finite dimensional graded $R$-modules, and hence all modules in 
$\Fdim R$. 

\subsection{}
We make $L$ a $\ZZ$-graded algebra by defining $\deg x_i=1$ and $\deg x_i^*=-1$ for all $i$.
The canonical map $R \to L$ is a homomorphism of graded rings and is injective. It is well-known, and not hard, to show that if $r>0$, then $L_r=x_0^rL_0$ and $L_{-r}=L_0(x_0^*)^r$ (see, e.g., \cite[Sect. 1.6]{Cuntz}). 
It follows from this that $L$ is strongly graded, i.e., $L_jL_{-j}=L_0$ for all integers $j$, and therefore
$$
\Gr L \equiv \Mod L_0
$$
where the functor giving the equivalence sends a graded module $M$ to its degree-zero component $M_0$. 

\begin{prop}
[Cuntz]
 \cite[Prop. 1.4]{Cuntz}
 $L_0 \cong S$. 
\end{prop}

\subsection{}
For $r \ge 1$, let $X_r$ be the set of words of length $r$ in the letters $x_i$, 
$0 \le i \le n$, and $X_\infty$ the union of all 
$X_r$, $r \ge 0$. 
Cuntz \cite[Lem. 1.3]{Cuntz}  shows that $L={\rm span}\{w^*w' \; | \; w,w' \in X_\infty\}$. 

\subsection{}
The next result was proved in \cite[Prop. 2.1]{A} but its utility is such that it seems useful to give a 
more direct proof.

\begin{prop}
Let $R=k\langle x_0,\ldots,x_n\rangle$ and define $L$ as above.
The ring $L$ is flat as a right $R$-module.
\end{prop}
\begin{pf}
We will show $L$ is an ascending union of finitely generated free right $R$-modules. Let
\begin{equation}
\label{defn.F_r}
F_r=\sum_{w \in X_r} w^* R.
\end{equation}
Suppose
$$
\sum_{w \in X_r} w^* r_w =0
$$
for some elements $r_w \in R$. Let $z \in X_r$. Then $zw^*=\d_{w,z}$ so $r_z=0$. It follows that all the $r_w$s are zero and the sum in  (\ref{defn.F_r}) is therefore a direct sum. Because $ww^*=1$, each $w^*R$ is a free $R$-module. Hence $F_r$ is free. 

Since $w^*=\sum_{i \in I} w^*x_i^* x_i$,  $F_r \subset F_{r+1}$. Since $L$ is spanned by elements $w^*w'$,  $L$ is the ascending union of the $F_r$s  and therefore flat. 
\end{pf}

A version of the following result for finitely presented not-necessarily-graded modules is 
given in \cite[Thm. 5.1]{A}.  Our proof differs in spirit from that in \cite{A}. 

\begin{prop}
\label{prop.L.tors}
Let $R$ and  $L$ be as above and $M \in \Gr R$. Then $L \otimes_R M=0$ if and only if $M \in \Fdim R$. 
\end{prop}
\begin{pf}
We observed in section \ref{sect.tors1} that $L \otimes_R M=0$ if $M \in \Fdim R$.

To prove the converse suppose $L \otimes_R M=0$.  First we will show $M$ is finite dimensional under the additional hypothesis that it is finitely presented. By Proposition \ref{prop.gr.coh}, 
$M_{\ge i}$ is a free $R$-module for $i \gg 0$. If $M_{\ge i}$ is a {\it non-zero} free module, then $L \otimes_R M$ would contain a non-zero free $L$-module. This does not happen because $L \otimes_R M =0$ so we deduce that $M_{\ge i}=0$ for $i \gg 0$. Since $M$ is finitely generated it is therefore finite dimensional.

Now let $M$ be an arbitrary graded $R$-module such that $L \otimes_R M=0$. To prove the proposition it suffices to show that $\dim_k Rm < \infty$ for all homogeneous $m \in M$.  Let $m \in M$ be a homogeneous element. 
Since $R$ is coherent $M$ is a direct limit of finitely presented graded modules, say $M=\liminj M_\lambda$ where each $M_\lambda$ is finitely presented. Let $\theta_\lambda:M_\lambda \to M$ be the canonical map and let $m_\lambda \in M_\lambda$ be
a homogeneous element such that $\theta_\lambda(m_\lambda)=m$.  Using $\theta_\lambda$, there is a map $L \otimes_R Rm_\lambda \to L \otimes_R Rm$; but $L \otimes_R Rm=0$
so the image of $1 \otimes m_\lambda$ in $L \otimes_R M_\nu$ is zero for some $\nu \gg i$. 
Let $m_\nu$ be the image of $m_\lambda$ in $M_\nu$. Then $L \otimes_R Rm_\nu=0$. Since $R$ is coherent and 
 $Rm_\nu$ is a finitely generated submodule of $M_\nu$, $Rm_\nu$ is finitely presented. The previous paragraph allows us to conclude that $\dim_k Rm_\nu < \infty$. But $Rm$ is the image of $Rm_\nu$ in $M$ so $\dim_k Rm<
 \infty$.  
\end{pf}

A version of the next result for finitely presented not-necessarily-graded modules is 
given in \cite[Thm. 5.1]{A}.  As with the previous result, the ideas in our proof 
 differ from those in \cite{A}.

\begin{thm}
\label{thm.Leavitt}
Let $\pi^*:\Gr R \longrightarrow \QGr R$ be the quotient functor and 
let $i^*=L \otimes_R - :\Gr R  \longrightarrow \Gr L$. Then 
$$
\QGr R \equiv \Gr L
$$
via a functor $\a^*:\QGr R \to \Gr L$ such that $\a^*\pi^* = i^*$.
\end{thm}
\begin{pf}
We already know $i^*$ is exact and vanishes on $\Fdim R$ so, by the universal property of  $\QGr R$, there is a unique functor $\a^*:\QGr R \to \Gr L$ such that $\a^*\pi^* = i^*$;
furthermore, $\a^*$ is exact.

The forgetful functor $i_*:\Gr L \to \Gr R$ is exact and right adjoint to $i^*$. We will show that $\pi^*i_*$ is quasi-inverse to $\a^*$. A diagram will help us keep track of the data:
$$
   \UseComputerModernTips
\xymatrix{
\Gr R \ar[r]^{\pi^*}   \ar[d]_{i^*} & \QGr R \ar[dl]^{\a^*}
\\
\Gr L \ar@/^2pc/[u]^{i_*}
}
$$
Since $R \to L$ is a universal localization it is an epimorphism in the category of rings. The multiplication map $L \otimes_R L \to L$ is therefore an isomorphism of $L$-bimodules. Thus, if $N \in \Gr L$, then
$$
i^*i_* N=L \otimes_R N = L \otimes_R(L \otimes_L N) \cong N.
$$
Therefore $\a^*(\pi^*i_*)=i^*i_* \cong \id_{\Gr R}$.

Let $M \in \Gr R$ and consider the exact sequence
 \begin{equation}
 \label{eq.4.term}
 0 \to \Tor_1^R(L/R,M) \to M \stackrel{f}{\longrightarrow} L \otimes_R M = i_*i^*M \to ( L/R) \otimes_R M \to 0
 \end{equation}
 where $f(m)=1 \otimes m$. Since $L \otimes_R L \cong L$, $i^*(f)$ is an isomorphism.
But $i^*$ is exact   so it   vanishes on $\Tor_1^R(L/R,M)$ and $( L/R) \otimes_R M$. 
Thus, by Proposition \ref{prop.L.tors}, both these modules are in $\Fdim R$. Therefore $\pi^*$ vanishes on them. Hence $\pi^*(f)$ is an isomorphism. In other words, the natural transformation $\pi^* \to \pi^*i_*i^*$ is an isomorphism.

In particular, $\pi^* \cong \pi^*i_*\a^*\pi^*$. By the universal property of $\QGr R$, there is a unique functor $\b^*$ such that the diagram
$$
   \UseComputerModernTips
\xymatrix{
\Gr R \ar[r]^{\pi^*}   \ar[d]_{\pi^*} & \QGr R \ar[dl]^{\b^*}
\\
\QGr R  
}
$$ 
commutes, i.e., $\pi^* = \b^*\pi^*$. That $\b^*$ is, of course, $\id_{\QGr R}$. But $\pi^*\cong ( \pi^*i_*\a^*)\pi^*$
so we conclude that $\pi^*i_*\a^* \cong \id_{\QGr R}$.  This completes the proof that $\a^*$ and $\pi^*i_*$
are mutually quasi-inverse. 
\end{pf}

\section{Remarks}
 
 \subsection{}
Most non-commutative projective algebraic geometry to date involves non-commutative rings that are noetherian. See, for example, Artin and Zhang's paper \cite{AZ} and the survey article of Stafford-Van den Bergh \cite{SV}. Two notable exceptions are (1) the rings $A:=k\langle x_1,\ldots, x_n \rangle/(f)$ where 
$f$ is a homogeneous quadratic  element of rank $n\ge 3$ (\cite{M0}, \cite{Pion1}) 
and (2) the non-commutative homogeneous coordinate rings appearing
 in Polishchuk's work  (\cite{Pol3} and \cite{Pol4}) on non-commutative elliptic curves 
 or, equivalently, non-commutative 2-tori endowed with a complex structure. 
The significance of the first is that $\sD^b(\QGr A)$ is equivalent to the bounded derived category of 
representations of the generalized Kronecker quiver (i.e., that with two vertices and $n$ parallel arrows from one to the other); $\Projnc A$ is viewed as a non-commutative analogue of the projective line.
Polishchuk's work provides a beautiful and deep connection between non-commutative geometry based 
on operator algebras and 
non-commutative projective algebraic geometry. 

The direct limit algebra $S$ in Theorem \ref{thm.3} provides a further link. When the base field is $\CC$ the norm-completion of $S$ belongs to an important class of $C^*$-algebras, the AF-algebras
(AF=approximately finite).  Under the philosophy that non-commutative $C^*$-algebras correspond to ``non-commutative topological spaces'' AF-algebras are often viewed as corresponding to 
0-dimensional spaces; see, for example, the paragraph at the foot of page 10 of Connes's 
book \cite{NCG}, although they also exhibit features of higher dimensional spaces. 
A prominent example is the AF-algebra associated to the space of Penrose tilings (see \cite{NCG} and \cite{Sm.y^2=0} for details).

\subsection{A homological remark}
A connected graded $k$-algebra $A$ is said to be {\sf Artin-Schelter regular} of dimension $d$ if
$\Ext^d_A(k,A) = k$ and $\Ext^i_A(k,A) =0$ when $i \ne d$. Many of the proofs in non-commutative projective algebraic geometry work only for Artin-Schelter regular rings of finite Gelfand-Kirillov
dimension.   

The next result shows that $TV$ is far from being Artin-Schelter regular when $\dim V \ge 2$.

\begin{lem}
Let $R$ be the free algebra on $d$ variables. Then 
there is an exact sequence 
$$
0 \to R^{d^2-1} \to  \uExt^1_R(k,R) \to k(1)^d \to  0
$$
of graded right $R$-modules. 
\end{lem}
\begin{pf}
Applying $\uHom_R(-,R)$ to a minimal resolution $0 \to R(-1)^d \to R \to k \to 0$ of left $R$-modules 
produces the top row in the commutative diagram of exact sequences
\begin{equation}
\label{Ext.diag}
   \UseComputerModernTips
\xymatrix{
 0 \ar[r] & \uHom(R,R) \ar[r]   \ar[d]^{\cong} & \uHom(R(-1)^d,R) \ar[r]   \ar[d]^{\cong} & \uExt^1_R(k,R)  \ar[r]  \ar@{=}[d]& 0
 \\
  0 \ar[r] & R \ar[r]_f &R(1)^d \ar[r] & \uExt^1_R(k,R)  \ar[r] & 0
  }
\end{equation}
 in which $f(1)=(x_1,\ldots,x_d)$ and $x_1,\ldots,x_d$ is a basis for $V$. 
 
 Let $e_1:=(1,0,\ldots,0), e_2:=(0,1,0,\ldots,0), \ldots, e_d:=(0,\ldots,0,1)$ be the standard basis for $R^d$.
Define the left $R$-module homomorphism $h:R \to (R^d)^{\oplus d}$ by $h(1)=(e_1,\ldots,e_d)$. 
Let  $g:R^d \to R(1)^d$ be the unique left $R$-module homomorphism such that $g(e_i) =x_i$ for all $i$.
There is an exact sequence $0 \to R^d \stackrel{g}{\longrightarrow} R(1) \to k(1) \to 0$.  
Let $(g,\ldots,g): (R^d)^{\oplus d} \to (R(1)^d)^{\oplus d}$ be the left $R$-module homomorphism 
defined by $(g,\ldots,g)(u_1,\ldots,u_d)=\big(g(u_1),\ldots,g(u_d)\big)$ where $u_i \in R^d$. 
Then there is a commutative diagram
\begin{equation}
\label{Ext.diag2}
   \UseComputerModernTips
\xymatrix{
& 0 \ar[d] & 0 \ar[d] 
\\
& R \ar[d]_h \ar@{=}[r] & R \ar[d]^f & 
\\
 0 \ar[r] & (R^d)^{\oplus d} \ar[r]^<<<<<<{(g,\ldots,g)}   \ar[d]    &R(1)^{\oplus d} \ar[r]   \ar[d]  & k(1)^{\oplus d} \ar[r]  \ar@{=}[d]& 0
 \\
& R^{d^2-1}  \ar[d]  & \uExt^1_R(k,R)   \ar[d]  & k(1)^d
\\
& 0  & 0
  }
\end{equation}
in which the columns are exact. The result now follows by applying the Snake Lemma to this diagram. 
 \end{pf}

The bottom row of (\ref{Ext.diag}) yields an exact sequence $0 \to \cO \to \cO(1)^d \to \cE \to 0$
in $\coh\XX^{d-1}$ (to be precise, since we started with left $R$-modules we should be replace the 
bottom row of (\ref{Ext.diag}) with the analogous exact sequence of left $R$-modules). By the lemma,
$\cE \cong \cO^{d^2-1}$). 
 
 \subsection{$\XX^n$ has only the trivial closed subspaces}
 
 There is a notion of closed subspace in non-commutative algebraic geometry \cite[Sect. 3.3]{VdB}. Rosenberg \cite[Prop. 6.4.1, p.127]{R} proved that closed subspaces of an affine nc-space 
 are in natural bijection with the two-sided ideals in a coordinate ring for it. The only two-sided ideals in 
  $S$ are the zero ideal and $S$ itself so the only closed subspaces of $\XX^n$ are the empty set and $\XX^n$ itself. 
 
 This is  a surprise because the free algebra contains a wealth of two-sided ideals. For example, 
the polynomial ring on $n+1$ variables is a quotient of the free algebra $k\langle x_0,\ldots,x_n\rangle$
so $\Qcoh \PP^n$ is a full subcategory of $\Qcoh \XX^n$ but coherent $\cO_{\PP^n}$-modules are not 
finitely presented as objects in $\Qcoh \XX^n$.
 
 \subsection{}
In \cite{Sm.kQ} we extend the ideas and results in this paper to path algebras of quivers:
the free algebra is replaced by a path algebra $kQ$ and  
the category of ``quasi-coherent sheaves'' on $\Projnc (kQ)$ is equivalent to the category of modules over a direct limit of semisimple $k$-algebras, namely $\liminj \End_{kI} (kQ_1)^{\otimes n}$ where $I$ is the set of vertices and $kQ_1$  the linear span of the arrows.

\end{document}